\documentclass[12pt]{article}
\usepackage{amsfonts}
\usepackage{amssymb,amsmath,amsthm,latexsym}
\usepackage[numbers,sort&compress]{natbib}
\usepackage{amsmath}
\usepackage{indentfirst}
\usepackage{adjustbox}

\usepackage[figuresright]{rotating}

\allowdisplaybreaks[4]
\newtheorem{thm}{Theorem}[section]
\newtheorem{prop}[thm]{Proposition}
\newtheorem{lem}[thm]{Lemma}
\newtheorem{cor}[thm]{Corollary}
\newtheorem{exam}{Example}
\newtheorem{defi}[thm]{Definition}
\newtheorem{rem}[thm]{Remark}

\newcommand{\pf}{{\bf Proof. \ }}

\numberwithin{equation}{section}

\begin{document}

\title{EAQEC codes from $s$-Galois hulls decomposition of linear codes}
\author{Hui  Li$^{1}$, Xiusheng Liu$^{2}$\footnote{Corresponding author.}
\setcounter{footnote}{-1}
\footnote{E-Mail addresses:
 christinehl405@163.com (H. Li),
 lxs6682@163.com (X. Liu)}}
\date{\small
1.School of Mathematics and Physics, Hubei Polytechnic University,\\ Huangshi, Hubei 435003, China \\
2. School of Science and Technology, College of Arts and Science of Hubei Normal University, Huangshi, Hubei 435109, China}
\maketitle


\begin{abstract} Entanglement-assisted quantum error-correcting (EAQEC)  codes are a subclass of quantum error-correcting codes
which use entanglement as a resource. These codes can provide error correction capability higher than the codes
derived from the traditional stabilizer formalism. In this paper, we first give a  $s$-Galois hulls   decomposition of linear codes over finite fields. By means of this decomposition, we provide a  method to construct EAQEC codes.  Using this method, we construct  new EAQEC  codes.

\end{abstract}


\bf Key Words\rm :  EAQEC ME code; $s$-Galois hull decomposition of a linear code;  matrix-product code

\bf MSC\rm:  81P70, 81P40, 94B15, 94B27.

\section{Introduction}
Nowadays quantum technologies become crucial to develop different areas of real world-life (see \cite{Ding,Dun,Sag}). So, quantum error-correcting codes are a necessary  tool in quantum computation and communication to detect and correct the quantum errors while quantum information is transferred via quantum channel.  One of the most known and used methods to create quantum error-correcting codes from classical block codes is the CSS method \cite{Calder1}. Unfortunately, it requires (Euclidean or Hermitian) duality containing to one of the classical codes used. One way to overcome this constraint is via entanglement. It is also possible to show that entanglement also improves the error-correction capability of quantum codes. These codes are called EAQEC codes. The first proposals of EAQEC codes were done by Bowen \cite{Bowen01}. In the following, Brun et al. \cite{Brun01} proposed to share entanglement between encoder and decoder to simplify the theory of quantum error correction and increase the communication capacity. With this new formalism, entanglement-assisted quantum stabilizer codes can be constructed from any classical linear code giving rise to EAQEC codes.  Fujiwara et al. \cite{Fuji} gave a general method for constructing entanglement-assisted quantum low-density parity check  codes.  Fan et al. provided a construction of EAQEC MDS
codes with a small number $c$ of pre-shared maximally entangled states in \cite{FanJ}. Li et al. in \cite{Li} constructed some EAQEC codes with good parameters by a decomposition of the defining set of a cyclic code. This method was generalized to apply  in negacyclic codes by Chen et al. \cite{Chen1}, where four classes of optimal EAQEC codes and two classes of EAQEC maximal-entanglement codes were obtained.  Lu et al. in \cite{Lu1} proposed a decomposition of the defining set of constacyclic codes and constructed four classes of EAQEC MDS codes from classical constacyclic MDS codes. Liu et al. in \cite{Liu} constructed six new classes of maximal-entanglement EAQEC codes by means of $s$-Galois dual codes.  Hu and Liu  in \cite{Hu} gave a new formula for calculating
the number $c$ of pre-shared maximally entangled states by using generator matrices of one code and  $s$-Galois dual code of the other code. By this formula, they constructed three classes of new EAQEC MDS codes.   Luo et al. in \cite{LM} gave  that EAQEC codes can be constructed using  linear codes, and  many authors have applied this result to construct EAQEC codes by non-binary linear codes (see \cite{Guenda,Hu1,Jin,Lu,Lu2,Luo,LiuH, Mu,Per,Pang,Wang,Wang1}).

In this paper, for any  linear code $C$  over finite fields $\mathbb{F}_{q}$,  we first present a $s$-Galois hull  decomposition  of the linear code $C$. Then we give a theorem of constructing   EAQEC  codes by the $s$-Galois hull decomposition of linear codes.

The paper is organized as follows. In Sect.$2$, we recall some basic knowledge on  linear codes and EAQEC codes. In Sect.$3$,  we  give  a $s$-Galois hull decomposition of the linear code $C$.  By means of this decomposition, we obtain a theorem of constructing  EAQEC  codes. Then  we construct  new EAQEC codes  in Sect.$4$.  Finally, a brief summary of this work is described in Sect. $5$.

\section{Preliminaries}
In order for the exposition in this paper to be self-contained, we introduce some basic concepts and results about linear codes and EAQEC codes, necessary for the development of this work. For more details, the interested reader is referred  to Refs. \cite{Mac, Liu1, Fuji, Brun01, Bowen01, Lai01,Lai02, Wilde}.

Throughout this paper, we denote the finite field of order $q=p^m$ by $\mathbb{F}_{q}$, where $p$ is
a prime and $m$ is a positive integer.  Let $\mathbb{F}_{q}^{*}$ be the multiplicative group of units of $\mathbb{F}_{q}$.
\subsection{Linear codes over $\mathbb{F}_q$}
For a positive integer $n$, let $\mathbb{F}_{q}^n=\{\mathbf{x}=(x_1,\cdots,x_n)\,|\,x_j\in \mathbb{F}_{q}\}$
which is an $n$ dimensional vector  space over $\mathbb{F}_{q}$.  A linear $[n,k]_{q}$ code $C$ over $\mathbb{F}_{q}$ is  an $k$ dimensional subspace of $\mathbb{F}_{q}^{n}$.  A generator matrix of an $[n, k]$ code $C$ is a $k \times n$ matrix whose rows are basis of $C$. The support of a vector $\mathbf{x}= (x_1, x_2,\ldots, x_n) \in \mathbb{F}_{q}^n$ is the subset $ \{j | x_j\neq 0\}$ of
$\{1, 2,\ldots, n\}$. The weight of a vector $\mathbf{x}=(x_1, x_2,\ldots, x_n) \in \mathbb{F}_{q}^n$ is the cardinality of the support of $\mathbf{x}$. A vector of $C$ is called a codeword of $C$. The minimum nonzero weight of all codewords in the linear $C$ is called
the minimum  Hamming weight $d_H(C)$ of $C$. An $[n, k]_q$ linear code with minimum weight $d$ is called an $[n, k, d]_q$ code.  One of the trade-offs between these parameters can be seen from the well-known Singleton bound:
$$d \leq n-k + 1.$$
An $[n, k, d]_q$ linear code $C$ achieving the Singleton bound is called a maximum distance separable (MDS) code.

For an $l\times n$ matrix $A =(A_{ij})_{l\times n}$ over $\mathbb{F}_{q}$, where $A_{ij}\in \mathbb{F}_{q}$. The transpose matrix of $A$,
denoted by $A^{T}$, is a matrix whose rows are the columns of $A$, i.e., $(A^T)_{ij}=A_{ji}$. Let $A^{(p^{m-s})} := (a_{ij}^{p^{m-s}})_{l\times n}$.
For a vector $\mathbf{a}=(a_1,a_2,\ldots,a_n)\in\mathbb{F}_{q} ^n$,
we have
$$\mathbf{a}^{p^{m-s}}=(a_1^{p^{m-s}},a_2^{p^{m-s}},\ldots,a_n^{p^{m-s}}).$$

Given two vectors $\mathbf{a}=(a_1,\ldots,a_n)$, and $\mathbf{b}=(b_1,\ldots,b_n)\in\mathbb{F}_{q}^n$, Fan and Zhang in \cite{Fan} introduced a kind of forms on $\mathbb{F}_{q}^n$
as follows {\rm(\cite[Definition 4.1]{Fan})}: for each integer $0\leq s< m $, define:
$$
[{\bf x},{\bf y}]_{s}=x_1y_1^{p^{s}}+\cdots+x_ny_n^{p^{s}},
\qquad\forall~ {\bf x},{\bf y}\in\mathbb{F}_{q}^n.
$$
We call $[{\bf x},{\bf y}]_{s}$ the $s$-Galois form on $\mathbb{F}_{q}^n$.
When $s=0$, it is just the usual Euclidean inner product.
When $m$ is even and $s=\frac{m}{2}$, it is  the Hermitian inner product.

Let $C$ be a linear code of length $n$  over $\mathbb{F}_{q}$. The $s$-Galois dual code of $C$
is defined as
$$C^{\bot_{s}}=\big\{{\bf x}\in\mathbb{F}_{q} ^n\,\big|\,
 [{\bf c},{\bf x}]_{s}=0,\, \forall~{\bf c}\in C\big\}
.$$
Then $C^{\bot_{0}}$ (simply, $C^{\perp}$) is just the Euclidean dual code of $C$,
and $C^{\bot_{\frac{e}{2}}}$ (simply, $C^{\perp_{H}}$)
is just the Hermitian dual code of $C$.

Obviously, $\dim_{\mathbb{F}_{q}}C+\dim_{\mathbb{F}_{q}}C^{\perp_s}=n$.

If $C \subset C^{\bot_{s}}$, then $C$ is said to be $s$-Galois self-orthogonal.
Moreover, $C$ is said to be $s$-Galois self-dual if $C= C^{\bot_{s}}$.


The following lemma can be found in \cite{Liu2}.
\begin{lem} \label{le:2.3} {\rm(\cite{Liu2}, Lemma $2.3$)}
Let $C$ be an $[n,k,d]$ linear code over $\mathbb{F}_{q}$ with a generator matrix $G$.
Then $C^{p^{m-s}}$ is an $[n,k,d]$ linear code  over $\mathbb{F}_{q}$
with a generator matrix $G^{(p^{m-s})}$;
and $(C^{p^{m-s}})^{\perp}=C^{\perp_{s}}$.
\end{lem}

The $s$-Galois hull of a linear code $C$ is defined as $\mathrm{Hull}_s(C) = C \cap C^{\perp_s}$.

A  linear code $C$ over $\mathbb{F}_{q}$ is called a
{\em linear complementary $s$-Galois dual code}
(abbreviated to {\em $s$-Galois LCD code})
if $C^{\perp_{s}}\cap C=\{\mathbf{0}\}$ (see \cite{Liu2}).

The following theorem shows a criteria of  $s$-Galois LCD codes which come from \cite{Liu2}.

\begin{lem} \label {le:2.4}{\rm(\cite{Liu2}, Theorem $2.4$)}
Let C be an $[n,k,d]_q$ linear code over $\mathbb{F} _{q}$ with generator matrix $G$.
Then $C$ is a $s$-Galois LCD code if and only if  $G(G^{(p^{m-s})})^{T}$ is nonsingular.
\end{lem}

\subsection{EAQEC codes over $\mathbb{F}_q$}
Let  $\mathbb{C}$ be the complex number field and  $V_n=\mathbb{C}^{q^{n}}=\mathbb{C}^{q}\otimes\cdots\otimes\mathbb{C}^{q}$ be $q^n$
dimensional Hilbert space. A quantum state is called a qudit which is a nonzero vector of $\mathbb{C}^{q}$.  An $n$-qudit is a nonzero vector in the $q^n$ dimensional Hilbert space $V_n$,  which can be expressed as
$$\{|\mathbf{v}\rangle=|v_1v_2\cdots v_n\rangle=|v_1\rangle \otimes|v_2\rangle \otimes\cdots\otimes| v_n\rangle:~\mathbf{v}=(v_1,v_2,\ldots,v_n) \in \mathbb{F}_{q}^{n}\}.$$


For $\mathbf{a}=(a_1,\ldots,a_n),\mathbf{b}=(b_1,\ldots,b_n),\mathbf{c}=(c_1,\ldots,c_n)\in \mathbb{F}_{q}^n$, define the unitary operators on $V_n$ as $X(\mathbf{a})|\mathbf{c}\rangle=|\mathbf{c}+\mathbf{a}\rangle$ and $Z(\mathbf{b})|\mathbf{c}\rangle=w^{tr[\mathbf{b},\mathbf{c}]}|\mathbf{c}\rangle$ where $w=\mathrm{exp}(2\pi i/p)$ is a primitive $p$-th root of unity and $\mathrm{tr}$ is the trace map from  $\mathbb{F}_{q}$ to  $\mathbb{F}_{p}$, i.e., $\mathrm{tr}(a)=a+a^{p}+a^{p^2}+\cdots+a^{p^{m-1}}$ for $a\in\mathbb{F}_{q}$.

Set $G_n=\{w^uX(\mathbf{a})Z(\mathbf{b}): \mathbf{a},\mathbf{b}\in \mathbb{F}_{q}^n,u\in\mathbb{F}_{p} \}$ is the error group associated with $E_n$. For any error $\mathbf{e} = w^uX(\mathbf{a})Z(\mathbf{b})$, its quantum weight is defined as $W_Q(\mathbf{e}) = |\{i|(a_i, b_i)\neq(0, 0)\}|$.

An  EAQEC code over $\mathbb{F}_{q}$, denoted by $[[n,k,d;c]]_q$, encodes $k$ information qudits into  $n$ channel qudits with the help of $c$ pairs of maximally-entangled states ($c$ ebits) and can correct up to $\lfloor \frac{d-1}{2} \rfloor$ errors, where $d$ is the minimum distance
of the code. 

If $c=0$, then the EAQEC code is a  quantum code. So,  EAQEC codes can be regarded as generalized quantum codes.



The following proposition is an important method  constructed EAQEC codes which can be found in \cite{Liu1}.

\begin{prop}\label{pro:2.4}{\rm(\cite{Liu1},Proposition 3.3)}
Let $C$ be an $[n, k, d]_q$ $s$-Galois LCD code. Then there exist an $[[n, k, d; n-k)]]_q$  EAQEC code.
\end{prop}




At the end of this section, we introduce another a definition on EAQEC codes.

\begin{defi}\label {de:2.2}
Let $C$ be an EAQEC code with parameters $[[n,k,d;c]]_{q}$. If  $c=n-k$, it is called an EAQEC maximal-entanglement (EAQEC ME) code.
\end{defi}

\section{$s$-Galois hulls decomposition of linear codes over $\mathbb{F}_q$}
In this section, we  provide a $s$-Galois hulls decomposition  of linear codes.

\begin{thm} \label{th:A} Let $C$ be an $[n,k,d]_q$ linear code over $\mathbb{F}_q$. Then $$C=\mathrm{Hull}_s(C)\oplus D,$$
where $D$ is a $s$-Galois LCD code over $\mathbb{F}_q$.
\end{thm}
\pf  Let $h = \mathrm{dim}_{\mathbb{F}_q}(\mathrm{Hull}_s(C))$ and $A = \{G_1, G_2,\ldots, G_h\}$ be a basis of $\mathrm{Hull}_s(C)$. Extend $A$ to be a basis $\{G_1, G_2,\ldots, G_h,G_{h+1},\ldots,G_k\}$ of $C$. Thus,
$$G_C=\begin{pmatrix}G_1\\G_2\\\vdots\\G_{h}\\G_{h+1}\\\vdots\\G_k\end{pmatrix},$$
is a generator matrix of $C$.

Take
$$G_D=\begin{pmatrix}G_{h+1}\\G_{h+2}\\\vdots\\G_k\end{pmatrix}.$$
Let $D$ be a linear code with the generator matrix $G_D$. Then $C=\mathrm{Hull}_s(C)\oplus D$.

The following we  prove  that the code $D$ is a $s$-Galois LCD  code. By Lemma \ref{le:2.4}, we only need to prove  $G_D(G_D^{(p^{m-s})})^{T}$ is nonsingular. We prove it by contradiction.

Obviously,
$$G_D(G_D^{(p^{m-s})})^T=\begin{pmatrix}G_{h+1}(G_{h+1}^{(p^{m-s})})^T&G_{h+1}(G_{h+2}^{(p^{m-s})})^T&\cdots&G_{h+1}(G_{k}^{(p^{m-s})})^T\\
G_{h+2}(G_{h+1}^{(p^{m-s})})^T&G_{h+2}(G_{h+2}^{(p^{m-s})})^T&\cdots&G_{h+2}(G_{k}^{(p^{m-s})})^T\\
\vdots&\vdots&\cdots&\vdots\\G_{k}(G_{h+1}^{(p^{m-s})})^T&G_{k}(G_{h+2}^{(p^{m-s})})^T&\cdots&G_{k}(G_{k}^{(p^{m-s})})^T\end{pmatrix}
=\begin{pmatrix}\mathbf{g}_{h+1}\\\vdots\\\mathbf{g}_{k}\end{pmatrix},$$
where $\mathbf{g}_j=(G_{j}(G_{h+1}^{(p^{m-s})})^T,G_{j}(G_{h+2}^{(p^{m-s})})^T,\cdots,G_{j}(G_{k}^{(p^{m-s})})^T)$ for $h+1\leq j\leq k$.

If  $G_D(G_D^{(p^{m-s})})^{T}$ is singular, then the vectors $\mathbf{g}_{h+1},\cdots,\mathbf{g}_{k}$ are linearly dependent. That is, for some $t$, $h+1\leq t \leq k$, there are $a_{h+1},\ldots,a_{t-1},a_{t+1},\ldots,a_k\in \mathbb{F}_q$ such that
$$\mathbf{g}_t=-a_{h+1}\mathbf{g}_{h+1}-\ldots-a_{t-1}\mathbf{g}_{t-1}-a_{t+1}\mathbf{g}_{t+1}+\ldots-a_k\mathbf{g}_{k}.$$
Without loss of generality, we assume that $t=k$. Then we have
$$\mathbf{g}_k=-a_{h+1}\mathbf{g}_{h+1}-\ldots-a_{k-1}\mathbf{g}_{k-1},$$
which  implies that
 $$\left\{\begin{aligned}&G_{k}(G_{h+1}^{(p^{m-s})})^T+a_{h+1}G_{h+1}(G_{h+1}^{(p^{m-s})})^T+a_{h+2}G_{h+2}(G_{h+1}^{(p^{m-s})})^T+\ldots-a_{k-1}G_{k-1}(G_{h+1}^{(p^{m-s})})^T=\mathbf{0},\\
&G_{k}(G_{h+2}^{(p^{m-s})})^T+a_{h+1}G_{h+1}(G_{h+2}^{(p^{m-s})})^T+a_{h+2}G_{h+2}(G_{h+2}^{(p^{m-s})})^T+\ldots+a_{k-1}G_{k-1}(G_{h+2}^{(p^{m-s})})^T=\mathbf{0},\\
&\cdots\cdots\cdots\cdots\cdots\cdots\cdots\cdots\cdots\cdots\cdots\cdots\cdots\cdots\cdots\cdots\cdots\cdots\cdots\cdots\cdots,\\
&G_{k}(G_{k-1}^{(p^{m-s})})^T+a_{h+1}G_{h+1}(G_{k-1}^{(p^{m-s})})^T+a_{h+2}(G_{h+2}G_{k-1}^{(p^{m-s})})^T+\ldots+a_{k-1}G_{k-1}(G_{k-1}^{(p^{m-s})})^T=\mathbf{0}.
\end{aligned}\right. $$
That is,
$$\left\{\begin{aligned}&G_{h+1}(G_{k}^{p^s})^T+a_{h+1}^{p^s}G_{h+1}(G_{h+1}^{p^s})^T+a_{h+2}^{p^s}G_{h+1}(G_{h+2}^{p^s})^T+\ldots-a_{k-1}^{p^s}G_{s+1}(G_{k-1}^{p^s})^T=\mathbf{0}^T,\\
&G_{h+2}(G_{k}^{p^s})^T+a_{h+1}^{p^s}G_{h+2}(G_{h+1}^{p^s})^T+a_{h+2}^{p^s}G_{h+2}(G_{h+2}^{p^s})^T+\ldots+a_{k-1}^{p^s}G_{h+2}(G_{k-1}^{p^s})^T=\mathbf{0}^T,\\
&\cdots\cdots\cdots\cdots\cdots\cdots\cdots\cdots\cdots\cdots\cdots\cdots\cdots\cdots\cdots\cdots\cdots\cdots\cdots\cdots\cdots,\\
&G_{k-1}(G_{k}^{p^s})^T+a_{h+1}^{p^s}G_{k-1}(G_{h+1}^{p^s})^T+a_{h+2}^{p^s}G_{k-1}(G_{h+2}^{p^s})^T+\ldots+a_{k-1}^{p^s}G_{k-1}(G_{k-1}^{p^s})^T=\mathbf{0}^T.
\end{aligned}\right. $$

It follows that
$$\begin{pmatrix}G_{h+1}\\\vdots\\G_{k}\end{pmatrix}((G_{k}+a_{h+1}G_{h+1}+a_{h+2}G_{h+2}+\ldots+a_{k-1}G_{k-1})^{p^s})^T=\mathbf{0}.$$
Clearly, $G_{k}+a_{h+1}G_{h+1}+a_{h+2}G_{h+2}-\ldots+a_{k-1}G_{k-1}\in C$, hence,
$$\begin{pmatrix}G_1\\\vdots\\G_h\\G_{h+1}\\\vdots\\G_{k}\end{pmatrix}((G_{k}+a_{h+1}G_{h+1}+a_{h+2}G_{h+2}+\ldots+a_{k-1}G_{k-1})^{p^s})^T=\mathbf{0}.$$
This means that $G_{k}+a_{h+1}G_{h+1}+a_{h+2}G_{h+2}+\ldots+a_{k-1}G_{k-1}\in \mathrm{Hull}_s(C)$.
Therefore, there are $a_{1},\ldots,a_h\in \mathbb{F}_q$ such that
$$G_{k}+a_{h+1}G_{h+1}+a_{h+2}G_{h+2}+\ldots+a_{k-1}G_{k-1}=a_{1}G_{1}+\ldots+a_hG_{h}.$$
That is,
$$G_{k}=a_{1}G_{1}+\ldots+a_hG_{h}-a_{h+1}G_{h+1}-a_{h+2}G_{h+2}-\ldots-a_{k-1}G_{k-1}.$$
This contradicts our assumption that $\{G_1, G_2,\ldots, G_h,G_{h+1},\ldots,G_k\}$ is a basis of $C$. Thus,  the matrix $G_D(G_D^{(p^{m-s})})^{T}$  is nonsingular, which is the expected result.
\qed

\section{An application EAQEC  codes}
In this section, we construct EAQEC  codes by using $s$-Galois hulls decomposition of linear codes.

Firstly, by means  of the Theorem \ref{th:A}, we give the following theorem constructed EAQEC codes.

\begin{thm} \label{th:B} Let $C$ be an $[n,k,d]_q$ linear code over $\mathbb{F}_q$. If $\mathrm{dim}_{\mathbb{F}_q}\mathrm{Hull}_s(C)=h$, then
there exist an $[[n, k-h, \geq d; n-k+h]]_q$  EAQEC code.
\end{thm}

\pf By Theorem \ref{th:A}, we have  $C=C_1\oplus C_2$, where $C_1=\mathrm{Hull}_s(C)$, and $C_2$ is a $s$-Galois LCD code.

Since $\mathrm{dim}_{\mathbb{F}_q}\mathrm{Hull}_s(C)=h$,  $C_2$ is a $s$-Galois LCD code  with parameters $[n,k-h,\geq d]_q$.
According to Proposition  \ref{pro:2.4},   there exist an $[[n, k-h, \geq d; n-k+h]]_q$  EAQEC code.
\qed
\subsection{Construction 1}
The following three lemmas can be found in \cite{Cao}.
\begin{lem} \label{le:4.1} {\rm( \cite{Cao}, Theorem III.1-III.3)} In $q=p^m$, let $p$ is  an odd prime.  Assume  $2s|m$. Let $n=\frac{r(q-1)}{p^s-1}$ for $1\leq r\leq p^s-1$.  Then, for $1\leq k \leq \lfloor\frac{p^s+n}{p^s+1}\rfloor$ and $0\leq h\leq k-1$,

$(1)$ there exists an $[n,k,n-k+1]_q$ code with $\mathrm{dim}_{\mathbb{F}_q}(Hull_s(C)) = h$.

$(2)$ there exists an $[n+1,k,n-k+2]_q$ code with $\mathrm{dim}_{\mathbb{F}_q}(Hull_s(C)) = h$.

$(3)$ there exists an $[n+2,k,n-k+3]_q$ code with $\mathrm{dim}_{\mathbb{F}_q}(Hull_s(C)) = h$.
\end{lem}

\begin{lem} \label{le:4.3} {\rm( \cite{Cao}, Theorem III.4-III.6)} In $q=p^m$, let $p$ is  an odd prime.  Assume  $2s|m, (q-1)|\mathrm{lcm}(x_1,x_2)$, and $\frac{q-1}{p^s-1}|x_1$ for two positive integers $x_1$ and $x_2$.  Let $n=\frac{r(q-1)}{\mathrm{gcd(x_1,q-1)}}$ for $1\leq r\leq\frac{q-1}{\mathrm{gcd(x_1,q-1)}}$.   Then, for $1\leq k \leq \lfloor\frac{p^s+n}{p^s+1}\rfloor$ and $0\leq h\leq k-1$,

$(1)$ there exists an $[n,k,n-k+1]_q$ code with $\mathrm{dim}_{\mathbb{F}_q}(Hull_s(C)) = h$.

$(2)$ there exists an $[n+1,k,n-k+2]_q$ code with $\mathrm{dim}_{\mathbb{F}_q}(Hull_s(C)) = h$.

$(3)$ there exists an $[n+2,k,n-k+2]_q$ code with $\mathrm{dim}_{\mathbb{F}_q}(Hull_s(C)) = h$.
\end{lem}

\begin{lem} \label{le:4.4} {\rm( \cite{Cao}, Theorem III.10-III.11)} In $q=p^m$, let $p$ is  an odd prime.  Assume  $2s|m$, and $a | s$. Let $n = tp^{aw}$ for each
$1 \leq t \leq p^{a}$ and each $1 \leq w \leq \frac{m}{a}-1$,  Then, for $1\leq k \leq \lfloor\frac{p^s+n-1}{p^s+1}\rfloor$ and $0\leq h\leq k-1$,

$(1)$ there exists an $[n,k,n-k+1]_q$ code with $\mathrm{dim}_{\mathbb{F}_q}(Hull_s(C)) = h$.

$(2)$ there exists an $[n+1,k,n-k+2]_q$ code with $\mathrm{dim}_{\mathbb{F}_q}(Hull_s(C)) = h$.
\end{lem}

\begin{cor} \label{co:4.1} In $q=p^m$, let $p$ is  an odd prime.   Assume  $2s|m$. Let $n=\frac{r(q-1)}{p^s-1}$ for $1\leq r\leq p^s-1$.  Then, for $1\leq k \leq \lfloor\frac{p^s+n}{p^s+1}\rfloor$ and $0\leq h\leq k-1$,

$(1)$ there exists an $[[n,k-h,\geq n-k+1,n-k+h]]_q$ EAQEC ME code.

$(2)$ there exists an $[[n+1,k-h,\geq n-k+2,n+1-k+h]]_q$ EAQEC ME code.

$(3)$ there exists an $[[n+2,k-h,\geq n-k+3,n+2-k+h]]_q$ EAQEC ME code.
\end{cor}

\pf $(1)$ By Lemma \ref{le:4.1}$(1)$, there exists an $[n,k,n-k+1]_q$ code with $\mathrm{dim}_{\mathbb{F}_q}(Hull_s(C)) = h$.
Again by Theorem  \ref{th:A}, there exists an $[n,k-h,\geq n-k+1]_q$ $s$-Galois LCD code.

According to Theorem \ref {th:B} and Definition \ref{de:2.2}, there exists an $[[n,k-h,\geq n-k+1,n-k+h]]_q$  EAQEC ME code.

$(2)$ and $(3)$ similar to the proof $(1)$.
\qed

\begin{exam} In Corollary \ref{co:4.1}, let $q = 3^4, s=1$, and $r=2,3$.  Then we obtain new EAQEC ME codes
in Tables $1$.
\end{exam}

\begin{table}
\caption{New EAQEC ME codes  from Corollary \ref{co:4.1}}
\begin{center}\begin{tabular}{cccc|cccc}
\hline
n& $k$& $h$& New EAQEC ME codes &n& $k$& $h$& New EAQEC ME codes\\
\hline
$312$&$14$&$1$&$[[312, 13, \geq 299; 299]]_{3^4}$ &$468$&$20$&$1$&$[[468, 19, \geq 449; 449]]_{3^4}$ \\
\hline
$312$&$15$&$1$&$[[312, 14, \geq 298; 298]]_{3^4}$& $468$&$21$&$1$&$[[468, 20, \geq 448; 448]]_{3^4}$ \\
\hline
\vdots&\vdots&\vdots&\vdots&\vdots&\vdots&\vdots&\vdots\\
\hline
$312$&$52$&$1$&$[[312,51, \geq 261; 261]]_{3^4}$ & $468$&$78$&$1$&$[[468, 77, \geq 391; 391]]_{3^4}$\\
\hline
$312$&$53$&$1$&$[[312,52, \geq 260; 260]]_{3^4}$ &$468$&$79$&$1$&$[[468, 78, \geq 390;390]]_{3^4}$ \\
\hline
$313$&$14$&$1$&$[[313, 13, \geq 300; 300]]_{3^4}$ &$469$&$20$&$1$&$[[469, 19, \geq 450; 450]]_{3^4}$ \\
\hline
$313$&$15$&$1$&$[[313, 14, \geq 299; 299]]_{3^4}$& $469$&$21$&$1$&$[[469, 20, \geq 449; 449]]_{3^4}$ \\
\hline
\vdots&\vdots&\vdots&\vdots&\vdots&\vdots&\vdots&\vdots\\
\hline
$313$&$52$&$1$&$[[313,51, \geq 262; 262]]_{3^4}$ & $469$&$78$&$1$&$[[469, 77, \geq 392; 392]]_{3^4}$\\
\hline
$313$&$53$&$1$&$[[313,52, \geq 261; 261]]_{3^4}$ &$469$&$79$&$1$&$[[469, 78, \geq 391;391]]_{3^4}$ \\
\hline
$314$&$14$&$1$&$[[314, 13, \geq 301; 301]]_{3^4}$ &$470$&$20$&$1$&$[[470, 19, \geq 451; 451]]_{3^4}$ \\
\hline
$314$&$15$&$1$&$[[314, 14, \geq 300; 300]]_{3^4}$& $470$&$21$&$1$&$[[470, 20, \geq 450; 450]]_{3^4}$ \\
\hline
\vdots&\vdots&\vdots&\vdots&\vdots&\vdots&\vdots&\vdots\\
\hline
$314$&$52$&$1$&$[[314,51, \geq 263; 263]]_{3^4}$ & $470$&$78$&$1$&$[[470, 77, \geq 393; 393]]_{3^4}$\\
\hline
$314$&$53$&$1$&$[[314,52, \geq 262; 262]]_{3^4}$ &$470$&$79$&$1$&$[[470, 78, \geq 392;392]]_{3^4}$ \\
\hline

\end{tabular}\end{center}
\end{table}

Similarly to Corollary \ref{co:4.1}, we get the following result by using Lemma \ref{le:4.3} and Theorem \ref {th:B}.
\begin{cor} \label{co:4.2} In $q=p^m$, let $p$ is  an odd prime.  Assume  $2s|m, (q-1)|\mathrm{lcm}(x_1,x_2)$, and $\frac{q-1}{p^s-1}|x_1$ for two positive integers $x_1$ and $x_2$.  Let $n=\frac{r(q-1)}{\mathrm{gcd(x_1,q-1)}}$ for $1\leq r\leq\frac{q-1}{\mathrm{gcd(x_1,q-1)}}$.   Then, for $1\leq k \leq \lfloor\frac{p^s+n}{p^s+1}\rfloor$ and $0\leq h\leq k-1$,

$(1)$ there exists an $[[n,k-h,\geq n-k+1,n-k+h]]_q$ EAQEC ME code.

$(2)$ there exists an $[[n+1,k-h,\geq n-k+2,n+1-k+h]]_q$ EAQEC ME code.

$(3)$ there exists an $[[n+2,k-h,\geq n-k+3,n+2-k+h]]_q$ EAQEC  ME code.
\end{cor}

\begin{exam} In Corollary \ref{co:4.2}, let $q = 3^6$, $s=1,x_1=364$,and $x_2=24$.  Then we obtain new EAQEC ME codes
in Tables $2$.
\end{exam}

\begin{table}
\caption{New EAQEC ME codes  from Corollary \ref{co:4.2}}
\begin{center}\begin{tabular}{cccc|cccc}
\hline
n& $k$& $h$& New EAQEC ME codes &n& $k$& $h$& New EAQEC ME codes\\
\hline
$91$&$6$&$1$&$[[91, 5, \geq86; 86]]_{3^6}$ &$182$&$9$&$1$&$[[182, 8, \geq 174; 174]]_{3^6}$ \\
\hline
$91$&$7$&$1$&$[[91, 6, \geq 85; 85]]_{3^6}$& $182$&$10$&$1$&$[[182, 9, \geq 173; 173]]_{3^6}$ \\
\hline
$91$&$8$&$1$&$[[91, 7, \geq 84; 84]]_{3^6}$& $182$&$11$&$1$&$[[182, 10, \geq 172;172]]_{3^6}$ \\
\hline
\vdots&\vdots&\vdots&\vdots&\vdots&\vdots&\vdots&\vdots\\
\hline
$91$&$23$&$1$&$[[91, 22, \geq 69; 69]]_{3^6}$ & $182$&$46$&$1$&$[[182, 45, \geq 137; 137]]_{3^6}$\\
\hline
$91$&$24$&$1$&$[[91,23, \geq 68;68]]_{3^6}$ &$182$&$47$&$1$&$[[182, 46, \geq 136;136]]_{3^6}$ \\
\hline
$92$&$6$&$1$&$[[92, 5, \geq87; 87]]_{3^6}$ &$183$&$9$&$1$&$[[183, 8, \geq 175; 175]]_{3^6}$ \\
\hline
$92$&$7$&$1$&$[[92, 6, \geq 86; 86]]_{3^6}$& $183$&$10$&$1$&$[[183, 9, \geq 174; 174]]_{3^6}$ \\
\hline
$92$&$8$&$1$&$[[92, 7, \geq 85; 85]]_{3^6}$& $183$&$11$&$1$&$[[183, 10, \geq 173;173]]_{3^6}$ \\
\hline
\vdots&\vdots&\vdots&\vdots&\vdots&\vdots&\vdots&\vdots\\
\hline
$92$&$23$&$1$&$[[92, 22, \geq 70; 70]]_{3^6}$ & $183$&$46$&$1$&$[[183, 45, \geq 137; 137]]_{3^6}$\\
\hline
$92$&$24$&$1$&$[[92,23, \geq 69;69]]_{3^6}$ &$183$&$47$&$1$&$[[183, 46, \geq 137;137]]_{3^6}$ \\
\hline
$93$&$6$&$1$&$[[93, 5, \geq88; 88]]_{3^6}$ &$184$&$9$&$1$&$[[184, 8, \geq 176; 176]]_{3^6}$ \\
\hline
$93$&$7$&$1$&$[[93, 6, \geq 87; 87]]_{3^6}$& $184$&$10$&$1$&$[[184, 9, \geq 175; 175]]_{3^6}$ \\
\hline
$93$&$8$&$1$&$[[93, 7, \geq 86; 86]]_{3^6}$& $184$&$11$&$1$&$[[184, 10, \geq 174;174]]_{3^6}$ \\
\hline
\vdots&\vdots&\vdots&\vdots&\vdots&\vdots&\vdots&\vdots\\
\hline
$93$&$23$&$1$&$[[93, 22, \geq 71; 71]]_{3^6}$ & $184$&$46$&$1$&$[[184, 45, \geq 139; 139]]_{3^6}$\\
\hline
$93$&$24$&$1$&$[[93,23, \geq 70;70]]_{3^6}$ &$184$&$47$&$1$&$[[184, 46, \geq 1368;138]]_{3^6}$ \\
\hline
\end{tabular}\end{center}
\end{table}

Finally, by Lemma \ref{le:4.4} and Theorem \ref {th:B}, we obtain the following two families of EAQEC ME codes.

\begin{cor} \label{co:4.3}In $q=p^m$, let $p$ is  an odd prime.  Assume  $2s|m$, and $a | s$. Let $n = rp^{aw}$ for each
$1 \leq r \leq p^{a}$ and each $1 \leq w \leq \frac{e}{a}-1$,  Then, for $1\leq k \leq \lfloor\frac{p^s+n-1}{p^s+1}\rfloor$ and $0\leq h\leq k-1$,

$(1)$ there exists an $[[n,k-h,\geq n-k+1,n-k+h]]_q$ EAQEC ME code.

$(2)$ there exists an $[[n+1,k-h,\geq n-k+2,n+1-k+h]]_q$ EAQEC ME code.
\end{cor}

\begin{exam} In Corollary \ref{co:4.3}, let $q = 7^6$ or $q = 13^4$ , $s=1,a=1,w=2$, and $r=4$. Then we obtain new EAQEC ME codes
in Tables $3$.
\end{exam}

\begin{table}
\caption{New EAQEC ME codes  from Corollary \ref{co:4.3}}
\begin{center}\begin{tabular}{cccc|cccc}
\hline
n& $k$& $h$& New EAQEC ME codes &n& $k$& $h$& New EAQEC ME codes\\
\hline
$196$&$4$&$1$&$[[196, 3, \geq193;183]]_{7^6}$ &$383$&$4$&$1$&$[[383, 3, \geq 380; 380]]_{13^4}$ \\
\hline
$196$&$5$&$1$&$[[196, 4, \geq 192; 192]]_{7^6}$& $383$&$5$&$1$&$[[383, 4, \geq 379; 379]]_{13^4}$ \\
\hline
$196$&$6$&$1$&$[[196, 5, \geq 191; 191]]_{7^6}$& $383$&$6$&$1$&$[[383, 5, \geq 378;378]]_{13^4}$ \\
\hline
\vdots&\vdots&\vdots&\vdots&\vdots&\vdots&\vdots&\vdots\\
\hline
$196$&$25$&$1$&$[[196, 24, \geq 172; 172]]_{7^6}$ & $383$&$25$&$1$&$[[383, 24, \geq359; 359]]_{13^4}$\\
\hline
$196$&$26$&$1$&$[[196,25, \geq 171;171]]_{7^6}$ &$383$&$26$&$1$&$[[383, 25, \geq 358; 358]]_{13^4}$ \\
\hline
$197$&$4$&$1$&$[[197, 3, \geq194; 194]]_{7^6}$ &$384$&$4$&$1$&$[[384, 3, \geq 381; 381]]_{13^4}$ \\
\hline
$197$&$5$&$1$&$[[197, 4, \geq 193; 193]]_{7^6}$& $384$&$5$&$1$&$[[384, 4, \geq 380; 380]]_{13^4}$ \\
\hline
$197$&$6$&$1$&$[[197, 5, \geq 192; 192]]_{7^6}$& $384$&$6$&$1$&$[[384, 5, \geq379;379]]_{13^4}$ \\
\hline
\vdots&\vdots&\vdots&\vdots&\vdots&\vdots&\vdots&\vdots\\
\hline
$197$&$25$&$1$&$[[197, 24, \geq 173; 173]]_{7^6}$ & $384$&$25$&$1$&$[[384, 24, \geq 360; 360]]_{13^4}$\\
\hline
$197$&$26$&$1$&$[[197,25, \geq 172;172]]_{7^6}$ &$384$&$26$&$1$&$[[384, 25, \geq359;359]]_{13^4}$ \\
\hline
\end{tabular}\end{center}
\end{table}

\subsection{Construction 2}
In this subsection, we  construct EAQEC codes  by using  matrix-product codes.

We first recall some basic concepts and results about the matrix-product codes.

Let $A=(a_{ij})_{l\times t}$ be an $l\times t$ matrix over $\mathbb{F}_{q}$ with $l \leq t$. Let $\mathcal{C}_1,\cdots,\mathcal{C}_l$ be codes of length $n$ over $\mathbb{F}_{q}$. Denote by $\mathrm{diag}(\lambda_1,\lambda_2,\ldots,\lambda_l)$,  the $l\times l$ diagonal matrix whose diagonal entries are $\lambda_1,\lambda_2,\ldots,\lambda_l$, where $\lambda_1,\lambda_2,\ldots,\lambda_l\in\mathbb{F}_{q}$.

The matrix product codes $\mathcal{C}=[\mathcal{C}_1,\ldots,\mathcal{C}_l]A$ is the set of  all matrix product $[\mathbf{c}_1,\cdots,\mathbf{c}_l]A$, where $\mathbf{c}_i\in \mathcal{C}_i$ are $n\times1$ column vectors for $1\leq i\leq l$.    If $\mathcal{C}_1,\ldots,\mathcal{C}_l$ are all linear codes with generator matrices $G_1,\ldots,G_l$, respectively, then  $[\mathcal{C}_1,\ldots,\mathcal{C}_l]A$ is a linear code generated by the following  matrix
$$G=\begin{pmatrix}a_{11}G_1&a_{12}G_1&\cdots&a_{1t}G_1\\ a_{21}G_2& a_{22}G_2 & \cdots&a_{2t}G_2 \\ \vdots &\vdots&\cdots& \vdots\\a_{l1}G_l & a_{l2}G_l& \cdots & a_{lt}G_l\end{pmatrix}.$$

Let us denote by $\mathbf{R}_i =(a_{i,1}, \ldots , a_{i,t})$ the element of $\mathbb{F}_q^t$ consisting of the $i$-th row of the matrix  $A$, for $i = 1, \ldots , l$.   We denote by $U_{A}(k)$ the code generated by $\langle\mathbf{ R}_1, \ldots , \mathbf{R}_k\rangle$ in  $\mathbb{F}_q^t$, where $1\leq k\leq l$.

Let $A=(a_{i j})_{l \times t}$ be a matrix over  $\mathbb{F}_q$. If the rows of $A$ are linearly independent, then we say that $A$ is a full-row-rank (FRR) matrix.

In the following, we list two useful results on matrix-product codes, which can be found in \cite{BN01,Ozb}.
\begin{lem}\label{lem:Bn1} {\rm(\cite{Ozb}, Theorem)}
Assume the notations are given as above. Let $\mathcal{C}_i$ be an $[n,k_i,d_i]_{q}$ linear code for $1\leq i\leq l$ and $A=(a_{ij})_{l\times t}$ be an FRR matrix. Let $\mathcal{C}=[\mathcal{C}_1,\cdots,\mathcal{C}_l]A$. Then $\mathcal{C}$ is an $[nt,\sum_{i=1}^{l}k_i,d(\mathcal{C})]_{q}$ linear code. Moreover, we have
$$d(\mathcal{C})\geq \mathrm{min}\{d_1d(U_A(1)),d_2d(U_A(2)),\ldots,d_ld(U_A(l)\}.$$
\end{lem}

\begin{lem} \label{lem:Bn3}{\rm(\cite{BN01}, Proposition $6.2$)}
Let the notations be given as above. Let  $A$ be an  $l\times l$ nonsingular matrix and $C_1,\cdots,C_l$ be linear codes over $\mathbb{F}_q$. Then
$$([C_1,\ldots,C_l]A)^{\perp}=[C_1^{\perp},\ldots,C_l^{\perp}](A^{-1})^{T}.$$
Furthermore,
$$d(C^{\perp })\geq \mathrm{min}\{d_1^{\perp}d(U_{(A^{-1})^T}(1)),d_2^{\perp}d(U_{(A^{-1})^T}(2)),\ldots,d_l^{\perp}d(U_{(A^{-1})^T}(l)\},$$
where $d_1^{\perp}=d_H(C_1^{\perp}),\ldots,d_l^{\perp}=d_H(C_l^{\perp})$.
\end{lem}

\begin{thm}\label{th:5.1}
For $1\leq i\leq l$, suppose that  $C_i\subset\mathbb{F}_{q}^{n}$ is an $[n,k_i,d_i]_{q}$  code.  Let $C=[C_1,C_2,\ldots,C_l]A $ where $A$ is $l\times l$ nonsingular matrix.  If $A (A^{(p^{m-s})})^{T}=\mathrm{diag}(\lambda_1,\lambda_2,\ldots,\lambda_l)$  with $\lambda_1,\ldots,\lambda_l\in\mathbb{F}_{q}^{*}$,  then

$(1)$ $\mathrm{Hull}_s(C)=[\mathrm{Hull}_s(C_1),\ldots,\mathrm{Hull}_s(C_l)]A $.

$(2)$   there exist an $[[2n, k-a, d; 2n-k+a]]_q$  EAQEC code, where $k=\sum_{i=1}^lk_i$, $a=\sum_{i=1}^l\mathrm{dim}_{\mathbb{F}_q}(\mathrm{Hull}_s(C_i)$,  and $d\geq \mathrm{min}\{d_1d(U_A(1)),d_2d(U_A(2)),\ldots,d_ld(U_A(l)\}$.
\end{thm}
\pf $(1)$ Since $A (A^{(p^{m-s})})^{T}=\mathrm{diag}(\lambda_1,\lambda_2,\ldots,\lambda_l)$ with $\lambda_1,\ldots,\lambda_l\in\mathbb{F}_{q}^{*}$,  we know that  $A^{(p^{m-s})}$ is  nonsingular.  Let $\Lambda=A (A^{(p^{m-s})})^{T}$. Then $(A^{(p^{m-s})})^{-1})^{T}=\Lambda^{-1}A$.

By Lemmas \ref{le:2.3}  and \ref{lem:Bn3},  we have
\begin{eqnarray*}
  ([C_1,\ldots,C_l]A)^{\perp_s}&=& (([C_1,\ldots,C_l]A)^{p^{m-s}})^{\perp}\\
   &=& [(C_1^{p^{m-s}})^{\perp},\ldots,(C_l^{p^{m-s}})^{\perp}]((A^{(p^{m-s})})^{-1})^{T}\\
   &=& [(C_1^{\perp_s},\ldots,C_l^{\perp_s}](\Lambda^{-1}A).
\end{eqnarray*}
As $\Lambda^{-1}=\begin{pmatrix}\lambda_{1}^{-1}&0&\cdots&0\\ 0&\lambda_{2}^{-1} & \cdots&0\\ \vdots &\vdots&\cdots& \vdots\\0 & 0& \cdots & \lambda_{l}^{-1}\end{pmatrix}$ and $\lambda_j^{-1}C_{j}^{\perp_s}=C_j^{\perp_s}$ for $j=1,\ldots,l$,  we have
$$[C_1^{\perp_s},\ldots,C_l^{\perp_s}]\Lambda^{-1}=[\lambda_{1}^{-1}C_1^{\perp_s},\ldots,\lambda_{l}^{-1}C_l^{\perp_s}]=[C_1^{\perp_s},\ldots,C_l^{\perp_s}].$$
Therefore,
$$C^{\perp_s}=[C_1^{\perp_s},\ldots,C_l^{\perp_s}]\Lambda^{-1}A=[C_1^{\perp_s},\ldots,C_l^{\perp_s}]A.$$

For any $\mathbf{w}\in \mathrm{Hull}_s(C)$, there are $\mathbf{a}_j\in C_j$ and $\mathbf{b}_j\in C_j^{\perp_s}$ for $j=1,\ldots,l$,   such that $\mathbf{w}=[\mathbf{a}_1,\ldots,\mathbf{a}_l]A= [\mathbf{b}_1,\ldots,\mathbf{b}_l]A$. Since $A$ is nonsingular,  $[\mathbf{a}_1,\ldots,\mathbf{a}_l]=[\mathbf{b}_1,\ldots,\mathbf{b}_l]$, which implies that $\mathbf{a}_1=\mathbf{b}_1,\ldots,\mathbf{a}_l=\mathbf{b}_l$. Hence, $\mathbf{w}\in [\mathrm{Hull}_s(C_1),\ldots,\mathrm{Hull}_s(C_l)]A $. This means that  $\mathrm{Hull}_s(C)\subset[\mathrm{Hull}_s(C_1),\ldots,\mathrm{Hull}_s(C_l)]A$.

Conversely, let $\mathbf{u}\in[\mathrm{Hull}_s(C_1),\ldots,\mathrm{Hull}_s(C_l)]A$, then there exist $\mathbf{c}_j\in \mathrm{Hull}_s(C_j)$ for  $j=1,\ldots,l$ such that $\mathbf{u}=[\mathbf{c}_1,\ldots,\mathbf{c}_l]A$.  Obviously, $\mathbf{u}=[\mathbf{c}_1,\ldots,\mathbf{c}_l]A\in \mathrm{Hull}_s(C)$. Thus, $[\mathrm{Hull}_s(C_1),\ldots,\mathrm{Hull}_s(C_l)]A\subset\mathrm{Hull}_s(C)$.

Summarizing, we have shown that $\mathrm{Hull}_s(C)=[\mathrm{Hull}_s(C_1),\ldots,\mathrm{Hull}_s(C_l)]A $.

$(2)$ By Theorem \ref{th:B} and Lemma \ref{lem:Bn1}, the statement is obvious.
\qed

\begin{cor}\label{co:5.1}
We assume that $C_i\subset\mathbb{F}_{q}^{n}$ is an $[n,k_i,d_i]_{q}$  code for $i=1,2$. Let $k=k_1+k_2$ and $a=\mathrm{dim}_{\mathbb{F}_q}(\mathrm{Hull}_s(C_1))+\mathrm{dim}_{\mathbb{F}_q}(\mathrm{Hull}_s(C_2))$. If $q = p^m$ where $p$ is an odd prime and $m$ is an odd positive integer, then there exist an $[[2n, k-a, d; 2n-k+a]]_q$  EAQEC code, where $d\geq \mathrm{min}\{2d_1,d_2\}$.
\end{cor}

\pf Let $C=[C_1,C_2]A$ where $A=\begin{pmatrix}1&1\\1&-1\end{pmatrix}$.
Obviously, $A(A^{(p^{m-s})})^{T}=\mathrm{diag}(2,2)$. By Theorem \ref{th:5.1}$(2)$, we know that there exists an EAQEC code with parameters $[[2n,k-a,d; 2n-k+a]]_q$, where $d\geq \mathrm{min}\{2d_1,d_2\}$.
\qed

\begin{exam} In Corollary \ref{co:5.1}, taking some special values of $q$, we obtain some new
EAQEC codes in Table $4$,  where $\alpha$ is a primitive element of $\mathbb{F}_{q}$.
\end{exam}
\begin{table}
   \caption{EAQEC ME codes from cyclic and negacyclic codes over $\mathbb{F}_{q}$}
\begin{center}\begin{tabular}{ccccc}
   \hline
$q$&$n$&Cyclic codes $C_1$& Negacyclic codes $C_2$ &New EAQEC ME codes\\
\hline
$9$&$8$&$C_1=\langle\Pi_{i= 2}^5(x-\alpha^i)$&$C_2=\langle\Pi_{i=0}^3 (x-\alpha^{1+2i})$&$[[16,3,\geq5;13]_9$ \\
\hline
$13$&$12$&$C_1=\langle\Pi_{i=3}^7(x-\alpha^i)$&$C_2=\langle \Pi_{i=0}^5 (x-\alpha^{1+2i})$&$[[24,5,\geq7;19]_{25}$\\
\hline
$25$&$24$&$C_1=\langle\Pi_{i=6}^{18}(x-\alpha^i)$&$C_2=\langle \Pi_{i=0}^{11} (x-\alpha^{1+2i})$&$[[48,11,\geq13;37]_{25}$\\
\hline
$29$&$28$&$C_1=\langle\Pi_{i=7}^{20}(x-\alpha^i)$&$C_2=\langle \Pi_{i=0}^{13} (x-\alpha^{1+2i})$&$[[56,13,\geq15;43]_{29}$\\
\hline
$49$&$48$&$C_1=\langle\Pi_{i=12}^{24}(x-\alpha^i)$&$C_2=\langle\Pi_{i=0}^{23 }(x-\alpha^{1+2i})\rangle$&$[[96,23,\geq25;73]_{49}$ \\
\hline
$169$&$168$&$C_1=\langle\Pi_{i=42}^{127}(x-\alpha^i)$&$C_2=\langle\Pi_{i=0}^{83} (x-\alpha^{1+2i})\rangle$&$[[336,83,\geq85;243]_{169}$\\
\hline
\end{tabular}\end{center}
\end{table}

In order to construct some new  EAQEC codes by utilizing Gabidulin codes, we need the following concepts and results which come from \cite{Gab,Hab,Jung}.

Let the elements of $\mathbf{g}=(g_1,\ldots,g_m)$  be a basis of  $\mathbf{F}_{q^m}$ over $\mathbb{F}_q$. The Gabidulin code $G_k(\mathbf{g})$ of length $m$ and dimension $k$ is defined by the generator matrix
$$~~~~~~~~~~~~~~M_{k}(g_1,g_2,\ldots,g_m)=\begin{pmatrix}g_1 & g_2 & \ldots &g_m \\
     g_1^{q} & g_2^{q} & \ldots &g_m^{q}\\
     g_1^{q^2} & g_2^{q^2} & \ldots &g_m^{q^2}\\
    \vdots& \vdots & \ddots & \vdots \\
   g_{1}^{q^{k-1}} & g_2^{q^{k-1}} & \ldots &g_m^{q^{k-1}}\\
  \end{pmatrix}. ~~~~~~~~~~~~~~~~~~~~~~~$$
It is well-known that the Gabidulin code $G_k(\mathbf{g})$ of length $m$ and dimension $k$ is an MDS code with parameters $[m,k,m-k+1]_{q^m}$ (see\cite{Gab}).

Let $\mathbf{g}=(g_1,\ldots,g_m)$ and $\mathbf{h}=(h_1,\ldots,h_m)$ be two bases of
$\mathbf{F}_{q^m}$ over $\mathbb{F}_q$. Then $\mathbf{h}$ is said to be a dual (orthogonal) basis of
$\mathbf{g}$ if $\mathrm{Tr}(g_ih_j) = \delta_{i,j}$, where  $\delta.,.$ is the Kroneker delta function and $\mathrm{Tr} : \mathbf{F}_{q^m}\rightarrow \mathbb{F}_q$ is the usual trace function. When $\mathbf{g}=\mathbf{h}$, we call $\mathbf{g}$ a self-dual basis.

The following result on self-dual basis is well-known.
\begin{prop} \label{pr:1} {\rm(\cite{Jung}, Theorem~1)}
In  $\mathbb{F}_{q^{m}}$, a  self-dual basis  $\mathbf{g}$ exists if and only if $q$ is even or both $q$ and $m$ are odd.
\end{prop}

In \cite{Hab}, authors provided an expression for calculating the hull dimension of Gabidulin codes.

\begin{lem}\label{le:5.2} Let $\mathbf{g}$ be a self-dual basis of $\mathbf{F}_{q^m}$ over $\mathbb{F}_q$, and
$G_k(\mathbf{g})$ be a Gabidulin code of length $m$ and dimension $k$. Then
$$\mathrm{dim}(\mathrm{Hull}_s(G_k(\mathbf{g}))=\left\{\begin{aligned}
         &\mathrm{min}\{m-k,s\}~~ ~~~~~~~\mathrm{if}~0\leq s\leq k,\\
         &\mathrm{min}\{m-s,k\}~~ ~~~~~~~\mathrm{if}~k+1\leq s\leq m-1.
         \end{aligned}\right. $$
\end{lem}

\begin{cor}\label{co:5.2}
For $q=2^m$, we assume that  $C_i\subset\mathbb{F}_{q}^{n}$ be an $[n,k_i,d_i]_{q}$  code for $i=1,2$. Let $k=k_1+k_2$ and $a=\mathrm{dim}_{\mathbb{F}_q}(\mathrm{Hull}_s(C_1))+\mathrm{dim}_{\mathbb{F}_q}(\mathrm{Hull}_s(C_2))$. Then there exist an $[[2n, k-a, d; 2n-k+a]]_q$  EAQEC code, where $d\geq \mathrm{min}\{d_1,d_2\}$.
\end{cor}

\pf Let $C=[C_1,C_2]A$ where $A=\begin{pmatrix}1&0\\0&1\end{pmatrix}$.
Clearly,  $A(A^{(2^{m-s})})^{T}=\mathrm{diag}(1,1)$. By Theorem \ref{th:5.1}$(2)$, we know that there exists an EAQEC code with parameters $[[2n,k-a,d; 2n-k+a]]_q$, where $d\geq \mathrm{min}\{d_1,d_2\}$.
\qed

Combining Corollary \ref{co:5.1} and Lemma \ref{le:5.2}, we have immediately following theorem.
\begin{thm}\label{th:5.2} We assume that $q$  is a power of an odd prime, and $m$ is an odd positive integer.  Let $\mathbf{g}=(g_1,\ldots,g_m)$ be a self-dual basis of $\mathbf{F}_{q^m}$ over $\mathbb{F}_q$, and $G_{k_1}(\mathbf{g})$ and $G_{k_2}(\mathbf{g})$ be two  Gabidulin codes of length $m$. Then

$(1)$ for $0 \leq s \leq \mathrm{min}\{k_1,k_2\}$, there exists an $[[2m,k-a,\geq d;2m-k+a]]_{q^m}$ EAQEC code, where $k=k_1+k_2, a = \mathrm{min}\{m-k_1,s\}+\mathrm{min}\{m-k_2,s\}$, and $d=\mathrm{min}\{2(m-k_1+1),m-k_2+1\}$.

$(2)$ for $0 \leq s \leq k_1, k_2 +1 \leq s \leq m-1$, there exists an $[[2m,k-a,\geq d;2m-k+a]]_{q^m}$ EAQEC code, where $k=k_1+k_2, a = \mathrm{min}\{m-k_1,s\}+\mathrm{min}\{m-s,k_2\}$, and $d=\mathrm{min}\{2(m-k_1+1),m-k_2+1\}$.

$(3)$ for $\mathrm{max} \{k_1+1, k_2 +1\} \leq s \leq m-1$, there exists an $[[2m,k-a,\geq d; 2m-k+a]]_{q^m}$ EAQEC code, where $k=k_1+k_2, a = \mathrm{min}\{m-s,k_1\}+\mathrm{min}\{m-s,k_2\}$, and $d=\mathrm{min}\{2(m-k_1+1),m-k_2+1\}$.
\end{thm}

Combining Corollary \ref{co:5.2} and Lemma \ref{le:5.2}, we have immediately following theorem.
\begin{thm}\label{th:5.3}  Let $\mathbf{g}=(g_1,\ldots,g_m)$ be a self-dual basis of $\mathbf{F}_{2^m}$ over $\mathbb{F}_2$, and
$G_{k_1}(\mathbf{g})$ and $G_{k_2}(\mathbf{g})$ be two  Gabidulin codes of length $m$. Then

$(1)$ for $0 \leq s \leq \mathrm{min}\{k_1,k_2\}$, there exists an $[[2m,k-a,\geq d;2m-k+a]]_{2^m}$ EAQEC code, where $k=k_1+k_2, a = \mathrm{min}\{m-k_1,s\}+\mathrm{min}\{m-k_2,s\}$, and $d=\mathrm{min}\{m-k_1+1,m-k_2+1\}$.

$(2)$ for $0 \leq s \leq k_1, k_2 +1 \leq s \leq m-1$, there exists an $[[2m,k-a,\geq d;2m-k+a]]_{2^m}$ EAQEC code, where $k=k_1+k_2, a = \mathrm{min}\{m-k_1,s\}+\mathrm{min}\{m-s,k_2\}$, and $d=\mathrm{min}\{m-k_1+1,m-k_2+1\}$.

$(3)$ for $\mathrm{max} \{k_1+1, k_2 +1\} \leq s \leq m-1$, there exists an $[[2m,k-a,\geq d; 2m-k+a]]_{2^m}$ EAQEC code, where $k=k_1+k_2, a = \mathrm{min}\{m-s,k_1\}+\mathrm{min}\{m-s,k_2\}$, and $d=\mathrm{min}\{m-k_1+1,m-k_2+1\}$.
\end{thm}

\begin{rem} \label{rem:5.1} In \cite{Hab,Hu},  authors also constructed EAQEC codes using the $s$-Galois hull dimension of the Gabidulin code, but their the Gabidulin code of length $m$ has dimension $k$ upper bounded by  $\lfloor\frac{p^s+m}{p^s+1}\rfloor$ or $\lfloor\frac{p^s+m+1}{p^s+1}\rfloor$
where $q=p^h$ and $2s | h$. Obviously,  in our case, dimension $k$ of the Gabidulin code of length $m$ has no such restriction, the codes parameters in Theorems \ref{th:5.2} and \ref{th:5.3} can exceed this bound. Therefore, we have more flexible parameters for EAQEC codes. On the other hand, lengths of EAQEC codes came from \cite{Hab,Liu} were $m$.  However, the lengths of three classes of EAQEC  codes derived from Theorems \ref{th:5.2} and \ref{th:5.3} are $2m$.
\end{rem}

\begin{exam} In Theorems \ref{th:5.2} and \ref{th:5.3}, taking some special values of $q$, we obtain some new
EAQEC codes in Table $5$.
\end{exam}
\begin{sidewaystable}
   \caption{EAQEC  codes from  Theorems \ref{th:5.2} and \ref{th:5.3}}
\scriptsize
   \begin{tabular}{llllllll}
   \hline
$(q=p^m,s)$&$k_1$&$k_2$&$a$& EAQEC  codes from  Theorems \ref{th:5.2} and \ref{th:5.3} &EAQEC  codes from  Refs &$k$&$\lfloor\frac{p^s+m}{p^s+1}\rfloor$\\
\hline
$(11^5,1)$&$k$&$k$&$2$&$[[10,2k-2,\geq 6-k,12-2k]]_{11^5}$&$[[5,2,3;1]_{11^5}$~\cite{Hu}&$2\leq k\leq 4 $&1 \\
\hline
$(13^6,2)$&$k$&$k$&$4$&$[[12,2k-4,\geq 7-k,14-2k]]_{13^6}$&$[[6,2,4;2]_{13^6}$~\cite{Hu}&$2\leq k\leq 4 $&1 \\
\hline
$(17^8,2)$&$k$&$k$&$2$&$[[16,2k-4,\geq 9-k,18-2k]]_{17^8}$&$[[8,4,4;2]_{17^8}$~\cite{Hu}&$2\leq k\leq 6 $&1 \\
\hline
$(2^{100},2)$&$k$&$k$&$2$&$[[200,2k-4,\geq 101-k,202-2k]]_{2^{100}}$&$[[100,k-2,101-k;98-k]]_{2^{100}}$~\cite{Hab}&$2\leq k\leq 98 $&20 \\
\hline
$(2^{100},2)$&$k$&$k$&$200-k$&$[[200,3k-200,\geq 101-k,400-4k]]_{2^{100}}$&$[[100,2k-100,101-k;0]]_{2^{100}}$~\cite{Hab}&$ 98\leq k< 101 $&20 \\
\hline
$(3^{67},40)$&$k$&$k$&$2k$&$[[134,0,\geq 68-k,134]]_{3^{67}}$&$[[67,0,68-k;67-2k]]_{3^{67}}$~\cite{Hab}&$ 2\leq k\leq 27 $&1 \\
\hline
$(3^{67},40)$&$k$&$k$&$200-2k$&$[[134,2(k-27),\geq 68-k,184-2k]]_{3^{67}}$&$[[67,k-27,68-k;40-k]]_{3^{67}}$~\cite{Hab}&$ 27\leq k\leq 39 $&1 \\
\hline
\end{tabular}
\end{sidewaystable}

\section{Conclusions}
In this paper, we exhibit the $s$-Galois hulls decomposition of linear codes over finite field $\mathbb{F}_{q}$.  As an application, we construct EAQEC  codes, among which some have dimensions $k > \lfloor\frac{p^s+m}{p^s+1}\rfloor$, and are hence not included in the codes that appeared in \cite{Cao, Hab,Hu}.

\section*{Data availability}

Data sharing not applicable to this article as no datasets were generated or analysed
during the current study.

\section*{Declaration of competing interest}

The author declare that she has no known competing financial interests or personal relationships that could have
appeared to influence the work reported in this paper.

\section*{Acknowledgements}
This work was supported by Scientific Research
Foundation of Hubei Provincial Education Department of China. (Grant
No. Q20174503) and the National Science Foundation of Hubei
Polytechnic University of China (Grant No. 17xjz03A).

\end{document}